\documentclass{roadef_en}

\usepackage[backend=biber, maxbibnames=99]{biblatex} 
\usepackage{graphicx}
\usepackage{fancyhdr}
\usepackage{hyperref}

\fancypagestyle{titlelogos}{%
  \fancyhf{}
  \fancyhead[L]{%
    \includegraphics[height=1.2cm]{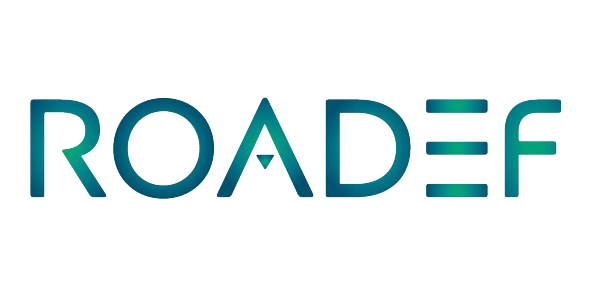}\hspace{0.5cm}%
    \includegraphics[height=1.2cm]{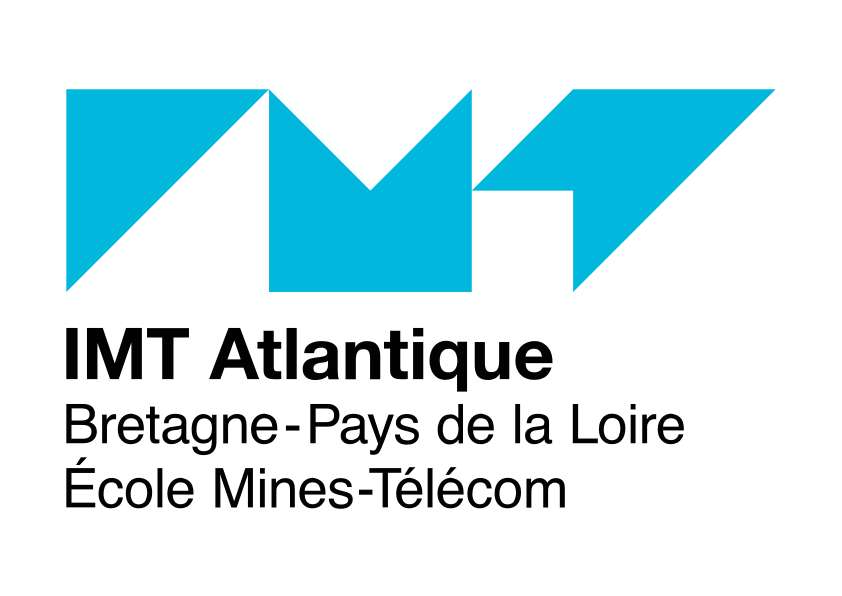}\hspace{0.5cm}%
    \includegraphics[height=1.2cm]{img/logo_insa_lyon.jpg}\hspace{0.5cm}%
    \includegraphics[height=1.2cm]{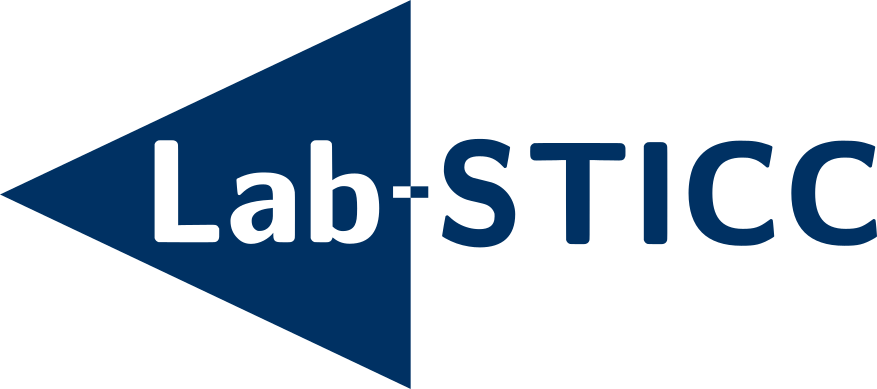}\hspace{0.5cm}%
    \includegraphics[height=1.2cm]{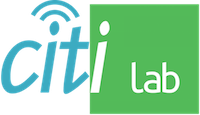}%
  }%
}

\begin{document}

\title{Branch \& Price \& Cut for the Time-Dependent Vehicle Routing Problem with Time Windows (TDVRPTW)}

\def\shorttitle{BPCA for the TDVRPTW}

\author{
    Florian Rascoussier\inst{1, 2}, 
    Romain Billot\inst{1},
    Lina Fahed\inst{1},
    Christine Solnon\inst{2}
}

\institute{
IMT Atlantique, Lab-STICC, UMR CNRS 6285, 29238 Brest, France  \\
\email{florian.rascoussier@imt-atlantique.fr}
\and
INSA Lyon, Inria, CITI, UR3720, 69621 Villeurbanne, France \\
\email{florian.rascoussier@insa-lyon.fr}
}

\maketitle
\thispagestyle{titlelogos}  


\section{Introduction}

In an urban context, travel times are highly dependent on the time of day, traffic conditions, and other external factors.
The \textbf{Time-Dependent Vehicle Routing Problem with Time Windows} (TDVRPTW) is an extension of the classical VRPTW, where a set of at most $k$ vehicles must perform routes to serve a set of customers while respecting strict time windows. Travel times depend on the departure time, and the objective is to minimize the sum of these travel times.

Column Generation (CG) allows decomposing this complex problem into two distinct and interdependent problems via Dantzig-Wolfe decomposition: given an incomplete set of candidate routes, the master problem is a linear Set Partitioning problem aiming to select at most $k$ routes such that each customer is visited by exactly one route, and the sub-problem (Pricing Problem) aims to calculate new routes to add to the set of candidate routes in order to improve the solution calculated by the master.

This method was used to solve the TDVRPTW in \cite{dabiaBranchPriceTimeDependent2013}, which uses CG within a Branch \& Price algorithm to solve the problem exactly. Their approach relies notably on an efficient DP algorithm to solve the pricing, using bidirectional search, relaxation of the elementarity constraint of customer visits, and the use of dominance criteria to limit the explosion of the number of states in memoization.

\section{Objectives}

Our work is part of the MAMUT project (\textit{Machine Learning and Matheuristics for Urban Transport}), \footnote{See the project page \url{https://anr.fr/Projet-ANR-22-CE22-0016}} which aims to develop hybrid solutions combining Operations Research (OR) techniques and machine learning to solve problems related to urban transport. Particular attention is paid to the explainability of algorithms, which restricts the use of heuristic or \textit{black-box} methods and orients the work towards interpretable approaches with guarantees on the quality of the solutions obtained.

In this context, the TDVRPTW problem seems particularly suitable due to its proximity to the real world and the various benchmarks available \cite{rifkiImpactSpatiotemporalGranularity2020}, but it remains relatively understudied. Currently, only\footnote{Post-comment: at the time of writing, the ...} the work of \cite{dabiaBranchPriceTimeDependent2013} constitutes the state of the art for the exact resolution of this problem, although their implementation is not accessible as open source. The first step is to reimplement this work with the aim of proposing a new version. Indeed, it seems that recent advances on the related problem of TDTSPTW, the single-vehicle version of the problem, could allow for significant improvements in results, particularly on the \textit{pricing} part. For this, we propose to use an exact and anytime extension of A* \cite{fontaineExactAnytimeApproach2023} and a \textit{Large Neighborhood Search} approach using dynamic programming to explore neighborhoods \cite{praletIteratedMaximumLarge2023a}, which have recently been successfully used to solve the TDTSPTW.

Finally, this work will explore the integration of data mining and machine learning techniques to guide column generation. These approaches aim to exploit knowledge from previously calculated routes to improve the efficiency of resolution algorithms. This last axis of improvement addresses the limitations of traditional OR methods, which ignore information available from past instances.

\section{Conclusion}

The TDVRPTW problem represents a realistic and complex extension of the VRPTW, taking into account the temporal constraints inherent in urban environments. This work is part of the MAMUT project, which aims to develop hybrid approaches combining advanced Operations Research techniques and machine learning tools to address the challenges posed by urban transport.

The main objective of this research is to improve the state of the art by bringing recent advances made on the TDTSPTW to the multi-vehicle context of the TDVRPTW. To this end, we propose to revisit the pioneering work of \cite{dabiaBranchPriceTimeDependent2013} by integrating modern methods, such as the Anytime A* algorithm \cite{fontaineExactAnytimeApproach2023}, and large neighborhood approaches based on dynamic programming \cite{praletIteratedMaximumLarge2023a}. These techniques, combined with the use of machine learning to effectively guide column generation, offer a unique opportunity to exploit both the strengths of explainable models and the richness of past data.

By exploring these complementary axes, this work aims to propose a high-performance and explainable \textit{solver}, capable of better meeting the requirements of real-world urban transport problems. These contributions will push the current limits of exact approaches and open new perspectives for the application of hybrid techniques in complex and dynamic contexts.

\printbibliography

\end{document}